\documentclass[twoside]{amsart}
\usepackage{latexsym}
\usepackage{amssymb,amsmath,amsopn}
\usepackage[dvips]{graphicx}   
\usepackage{color,epsfig}      
\begin{document}
{\Large Hamiltonian paths on directed grids}

\medskip Hujter, Mih\'aly \ and \ Kaszanyitzky, Andr\'as

Budapest University of Technology and Economics

\texttt{hujter@math.bme.hu; kaszi75@gmail.com}

\bigskip \noindent \textbf{Abstract.} Our studies are related to a special
class of FASS-curves, which can be described in a 
node-rewriting
Lindenmayer-system. These ortho-tile (or diagonal) type recursive curves
inducing Hamiltonian paths. We define a special directed graph on a
rectangular grid, and we enumerate all Hamiltonian paths on this graph. Our
formulas are strongly related to both the Fibonacci numbers and the domino
tilings of chessboards. The constructability of the regular $17$-gon with
straightedge and compass is also related.

\medskip \noindent \textbf{Introduction.} In 1877, G.\ Cantor proved for any
positive integer $d$ that, there exists a one-to-one point-to-point
correspondence between a unit line segment and the entire $d$%
-di\-men\-sio\-nal space. In other words, the infinite number of points in a
unit interval is the same~cardinality~as the infinite number of points in
any finite-dimensional manifold. In 1890, G.\ Peano constructed a continuous
mapping from the unit interval onto the unit square. This continuous curve
that passes through every point of the unit square was the first example of
space filling curves.

In 1891, D.\ Hilbert discovered another type of these recursive curves
[Sa94]. These FASS-curves (space-\textbf{F}illing, self-\textbf{A}voiding, 
\textbf{S}imple, self-\textbf{S}imilar) can be described in a node-rewriting
Lindenmayer-system [PL90,PLF91]. The \emph{ortho-tile} type of FASS-curves
can be represented only on a special directed grid graph. (See Figures 8.a
and 8.b.) Their approximations are Hamiltonian paths between diagonally
located points. 

\medskip \noindent \textbf{Definitions.} 
Given fixed positive integers $p,q$, 
by a \emph{directed grid graph with an odd-even direction}, or by $%
DGG_{p,q}$, in short, we mean a directed graph on the vertex set $%
\{1,2,...,p\} \times \{1,2,...,q\}$, as a subset of the $xy$ coordinate
plane, with all possible arcs of the form $(x;y)\rightarrow (x+1;y)$ if $y$
is odd, and arcs of the form $(x;y)\rightarrow (x-1;y)$ if $y$ is even, and
arcs of the form $(x;y)\rightarrow (x;y+1)$ if $x$ is odd, and arcs of the
form $(x;y)\rightarrow (x;y-1)$ if $x$ is even. This graph has $pq$ vertices
and $(p-1)q+p(q-1)=\allowbreak 2pq-p-q$ arcs. A \emph{Hamiltonian path} is
such a permutation $v_{1}\rightarrow v_{2}\rightarrow ...\rightarrow v_{pq}$
of all vertices for which $v_{n}\rightarrow v_{n+1}$ is an arc for each 
$n\in \{1,2,...,pq-1\}$. Let $h(p,q)$ denote the number of Hamiltonian paths
in $DGG_{p,q}$. Furthermore, for a positive integer $r$, let $h_{r}(p,q)$
denote the number of those Hamiltonian paths of the form $v_{1}\rightarrow
v_{2}\rightarrow ...\rightarrow v_{pq}$ for which $v_{1}=(1;1)$, $%
v_{2}=(2;1) $, ..., $v_{r}=(r;1)$ hold but $v_{r+1}\neq (r+1;1)$. For other
definitions and for the history of the Hamiltonian paths we the refer reader
to [Ba06] and [We01].

In Figure 1 we show $DGG_{5,4}$, and in Figures 2.a and 2.b we give two
Hamiltonian paths.

\medskip \noindent \textbf{Basic observations.} We can make the following
easy observations: For any positive integers $p,q$, we have $h(p,q)=h(q,p)$
and $h(p,1)=h(1,q)=1 $. If $p$ is odd, then $h(p,2)=h(2,p)=1$. If both $p$
and $q$ are even, then $h(p,q)=0$. If $pq$ is odd, then $h(p,q)>0$.

We can also observe that for any Hamiltonian path $v_{1}\rightarrow
v_{2}\rightarrow ...\rightarrow v_{pq}$ we have $v_{1}=(1;1)$. On the other
hand, if both $p$ and $q$ are odd, then $v_{pq}=(p;q)$. If $p$ is odd but $q$
is even, then $v_{pq}=(1;q)$. Finally, if $p$ is even but $q$ is odd, then $%
v_{pq}=(p;1)$. For example, in case of the Hamiltonian paths of Figures 2.a
and 2.b, the terminal points are in the top-left corner.

We can also make the following observations: If $r$ is even, then $%
h_{r}(p,q)=0$. Therefore%
\[
h(p,q)=h_{1}(p,q)+h_{3}(p,q)+h_{5}(p,q)+\cdots+h_{p}(p,q) 
\]%
if $p$ is odd, and 
\[
h(p,q)=h_{1}(p,q)+h_{3}(p,q)+h_{5}(p,q)+\cdots+h_{p-1}(p,q) 
\]%
if $p$ is even. If both $p$ and $q$ are odd and both are at least 3, then $%
h_{1}(p,q)>0$, $h_{3}(p,q)>0$, ..., $h_{p}(p,q)>0$. If $p$ is odd and $q\geq
2$, then $h_{p}(p,q)=h(p,q-1)$. If $q\geq 3$, then $h_{1}(3,q)=h(3,q-2)$.

\medskip \noindent \textbf{Fibonacci numbers.} Among the above observations
we find that $h(1,3)=h(3,1)=1$, $h(2,3)=h(3,2)=1 $, and for $n=3,4,5,...$,
we have 
\[
h(3,n)=h_{1}(3,n)+h_{3}(3,n)=h(3,n-2)+h(3,n-1). 
\]%
Therefore the numbers $h(1,3)=h(3,1)$, $h(2,3)=h(3,2)$, $h(3,3)$, ... are
exactly the well-known Fibonacci numbers: $F_{1}=1$, $F_{2}=1$, $F_{3}=2$, $%
F_{4}=3$, $F_{5}=5$, ..., $F_{n}=F_{n-1}+F_{n-2}$, ... . Up to now we gained
the following table of the values $h(p,q)$ for $\min \{p,q\} \leq 6$. 
\[
\begin{tabular}{|l|c|c|c|c|c|c|}
\hline
$\ $ & ${\small p=1}$ & ${\small p=2}$ & ${\small p=3}$ & ${\small p=4}$ & $%
{\small p=5}$ & ${\small p=6}$ \\ \hline
$h(p,1)=$ & 1 & 1 & 1 & 1 & 1 & 1 \\ \hline
$h(p,2)=$ & 1 & 0 & 1 & 0 & 1 & 0 \\ \hline
$h(p,3)=$ & 1 & 1 & 2 & 3 & 5 & 8 \\ \hline
$h(p,4)=$ & 1 & 0 & 3 & 0 &  & 0 \\ \hline
$h(p,5)=$ & 1 & 1 & 5 &  &  &  \\ \hline
$h(p,6)=$ & 1 & 0 & 8 & 0 &  & 0 \\ \hline
\end{tabular}%
\]%
In the rest of the present paper we focus on the cases where $\min \{p,q\}
\geq 4$.

\medskip \noindent \textbf{Definitions.} By a \emph{domino} we mean such a
rectangle whose corners are all vertices of our digraph, and in the
rectangle the length of the diagonals is exactly $\sqrt{5}$. We say that
some pairwise nonoverlapping dominoes form a \emph{domino tiling} of our
graph if their union is the entire grid, i.e. the rectangle with corners $%
(1;1)$, $(p;1)$, $(p;q)$, $(1;q)$. Clearly, the number of dominoes in a
domino tiling is $(p-1)(q-1)/2$. We say that a given Hamiltonian path $%
v_{1}\rightarrow v_{2}\rightarrow ...\rightarrow v_{pq}$ \emph{avoids} a
given domino if the shorter symmetry axis of the domino is not an arc of the
Hamiltonian path.

\medskip \noindent \textbf{Theorem 1. }\emph{There is a one-to-one
correspondence between Hamiltonian paths and the domino tilings such that
for each domino tiling the corresponding Hamiltonian path is the only one
which bisects no domino of the tiling.} \medskip

We will prove the theorem later. In Figures 3.a and 3.b we find the domino
tilings that correspond to the Hamiltonian paths of Figure 2.a and 2.b.,
respectively.

More than a half a century ago [Ka61] and [Te61] proved that the number of
different domino tilings is exactly%
\[
\prod_m\prod_k\left( 4\cos ^{2}\frac{m\pi }{p}+4\cos
^{2}\frac{k\pi }{q}\right) 
\]%
where the products are understood for all positive integers $m$ and $k\ $%
such that $2m<p$, $2k<q$. This formula is well-known as the \emph{Kasteleyn}
formula.

\medskip \noindent \textbf{Example.} We consider the case $p=q=5$. Now 
\[
4\cos ^{2}\frac{\pi }{5}=\frac{3+\sqrt{5}}{2} \  \  \  \  \  \  \  \  \  \  \ 4\cos
^{2}\frac{2\pi }{5}=\frac{3-\sqrt{5}}{2} 
\]%
\[
\Pi_m\Pi_k\left( 4\cos ^{2}\frac{m\pi }{p}
+4\cos^{2}\frac{k\pi }{q}\right) =\left( 3+\sqrt{5}\right) \cdot 3^{2}\cdot
\left(3-\sqrt{5}\right) =36 
\]%
Therefore, by Theorem 1 we gain that $h(5,5)=36$.

\medskip \noindent \textbf{Remark.} The above mentioned Kasteleyn formula
outputs $0$ if both $p $ and $q$ are even. If $\min \{p,q\}=1$, the meaning
of the formula is $1$. If \ $\min \{p,q\}=2$ and $\left \vert
p-q\right
\vert \in 1,3,5,...$, then the Kasteleyn formula outputs $1$.
Since $4\cos ^{2}\frac{\pi }{3}=$$1$, if \ $\min \{p,q\}=3$, then the
Kasteleyn formula and our Theorem 1 produces the following two nice formulas
for the Fibonacci numbers $F_{2n}$ and $F_{2n+1}$ for any positive integer $%
n $. 
\[
\prod_{k=1}^{n-1}\left( 1+4\cos ^{2}\frac{k\pi }{2n}\right)
=F_{2n}\  \  \  \  \  \  \  \prod_{k=1}^n\left( 1+4\cos ^{2}\frac{k\pi }{%
2n+1}\right) =F_{2n+1} 
\]%
For example for $n=8$ and for $x_{k}=1+4\cos ^{2}\frac{k\pi }{17}$, $%
k=1,2,...,8$, the latter formula produces $x_{1}x_{2}\cdots x_{8}=1597$.
However, by the famous result of Gauss on the constructability of the
regular $17$-gon with straightedge and compass, we can derive a nice formula
for each $x_{k}$ separately. Namely, we can start from the well-known
formula (see, e.g., [Gi07])%
\[
16\cos \frac{2\pi }{17}=\sqrt{17}-1+\sqrt{34-2\sqrt{17}}+2\sqrt{17+3\sqrt{17}%
-\sqrt{170+38\sqrt{17}}} 
\]%
and we can apply that 
\[
1+\cos \frac{2\pi }{17}=2\cos ^{2}\frac{\pi }{17} \  \  \  \  \  \  \  \  \  \  \  \  \
\  \ 1-\cos \frac{2\pi }{17}=2\sin ^{2}\frac{\pi }{17} 
\]%
This way we can express each $x_{k}$ as $\frac{1}{16}$ times an integer
coefficient polynomial of 
$v=\sqrt{34-2\sqrt{17}}$ and $w=\sqrt{17+3\sqrt{17}-\sqrt{170+38\sqrt{17}}}$ 
because 
$64\cos^2\frac{\pi }{17}\ =64+(2-v)v+4w$. 
For example $x_{1}=5+\frac{(2-v)v}{16}+\frac{w}{4}$. 

\medskip \noindent \textbf{Table of $h(p,q)$.} By the Kasteleyn formula and
by our Theorem 1 we can continue the above incomplete table as follows%
\[
\begin{tabular}{|l|c|c|c|c|c|c|}
\hline
$\ $ & ${\small p=4}$ & ${\small p=5}$ & ${\small p=6}$ & ${\small p=7}$ & $%
{\small p=8}$ & ${\small p=9}$ \\ \hline
$h(p,4)=$ & 0 & 11 & 0 & 41 & 0 & 153 \\ \hline
$h(p,5)=$ & 11 & 36 & 95 & 281 & 781 & 2245 \\ \hline
$h(p,6)=$ & 0 & 95 & 0 & 1183 & 0 & 14824 \\ \hline
$h(p,7)=$ & 41 & 281 & 1183 & 6728 & 31529 & 167089 \\ \hline
$h(p,8)=$ & 0 & 781 & 0 & 31529 & 0 & 1292697 \\ \hline
$h(p,9)=$ & 153 & 2245 & 14824 & 167089 & 1292697 & 12988816 \\ \hline
\end{tabular}%
\]

\medskip \noindent \textbf{Definition.} Given a domino $D$, there are six
vertices at the perimeter of the domino; the \emph{canonical numbering} of
these vertices is $D_{1},D_{2},...,D_{6}$ if all five $D_{j}\rightarrow
D_{j+1}$ are arcs in the graph, $j=1,2,...,5$. In Figure 4 we can see a
canonical numbering. Observe, that each domino has exactly one canonical
numbering. In this paper we always consider canonical numberings.

\medskip \noindent \textbf{The proof of Theorem 1.} The cases $\min \{p,q\}
\leq 2$ are obvious. In the rest of the proof we assume that $\min \{p,q\}
\geq 3$. The case when both $p$ and $q$ are even is also obvious. By
symmetry we may assume that $p$ is odd. We have at least one Hamiltonian
path and we have at least one domino tiling.

Let us consider a fixed Hamiltonian path $H$ and a fixed domino $D$.
According to the canonical numbering let the vertices around $D$ be $D_{i}$, 
$i=1,2,...,6$. Observe that if $D_{5}\rightarrow D_{2}$ is not an arc in $H$%
, then $D_{5}\rightarrow D_{6}$ must be an arc in $H$, and $D_{1}\rightarrow
D_{6}$ cannot be an arc in $H$. Given $H$ and given $k$ pairwise
nonoverlapping dominoes, we gain $2k$ arcs such that none of them can be in $%
H$. These $2k$ arcs are all distinct because the $D_{5}\rightarrow D_{2}$
arcs are all inside the pairwise nonoverlapping dominoes, and each $%
D_{1}\rightarrow D_{6}$ arc belongs to the domino situated in the right
angle determined by the two arcs started at $D_{1}$. In a domino tiling
there are $(p-1)(q-1)/2$ pairwise distinct dominoes. In summary, if $H$
avoids each domino in a fixed domino tiling, then there are $pq-1$ arcs in $%
H $ and there are $(p-1)(q-1)$ further arcs not in $H$. However 
\[
pq-1+(p-1)(q-1)=2pq-q-p 
\]%
This is exactly the total number of arcs in the graph. Therefore, starting
from the given domino tiling, some (but not necessary all) of the canonical
numbering $D_{j}\rightarrow D_{j+1}$ arcs of the domino's perimeter form the
Hamiltonian path. On the other hand, the $(p-1)(q-1)$ arcs missing from the
a Hamiltonian path determine two-by-one all dominoes of the domino tiling.

Now we make the above argument more explicit. As on a usual chessboard, we
say that a unit square of our grid is black if the sum of bottom-left corner
coordinates is even. (See Figure 5 as an illustration.) The other unit
squares are called white. Clearly, each domino $D$ consists of a black
square and of a white square; the corners of the black square are $D_{1}$, $%
D_{2}$, $D_{5}$, $D_{6}$ according to the canonical ordering, and the white
square's corners are $D_{2}$, $D_{3}$, $D_{4}$, $D_{5}$ . If the domino is
an element of the domino tiling, than the arc between the domino's squares
must not be in the corresponding Hamiltonian path. These $(p-1)(q-1)/2$
pairwise distinct arcs of $DGG_{p,q}$ are called \emph{domino-axis} arcs.
According to the canonical ordering of a domino $D$, the domino-axis arc is
the $D_{5}\rightarrow D_{2}$ arc. On the other hand, if a Hamiltonian path
avoids a domino $D$, then the $D_{1}\rightarrow D_{6}$ arc can not be in the
Hamiltonian path, either. We call these arcs \emph{domino-black-end} arcs.
Clearly, each domino-axis arc corresponds to exactly one domino in a domino
tiling, and each domino-black-end arc corresponds to exactly one domino in a
domino tiling, too. The total number of domino-axis arcs and
domino-black-end arcs is $(p-1)(q-1)$. The number of the remaining arcs is
exactly the same as the number of the arcs in a Hamiltonian path. Therefore
a domino tiling determines at most one Hamiltonian path which bisects no
domino of the tiling.

We can observe that in case of a domino $D$ of a domino tiling, the arcs $%
D_{1}\rightarrow D_{2}$ and $D_{5}\rightarrow D_{6}$ can be neither
domino-axis arcs nor domino-black-end arcs.\ Therefore both arcs must be in
the Hamiltonian path corresponding to the domino tiling.

Now let us consider such arcs of the grid graph which are at the perimeter
of the grid and which are at the perimeter of a white square at the same
time. One can easily compute that there are $(p-1)/2$ such arcs at the
bottom line, $\left( p-1\right) /2$ such arcs at the top line, and there are
other $q-1$ such arcs at the vertical sides of the grid. We call such arcs
as \emph{white perimeter} arcs. Observe that each white perimeter arc must
be in each Hamiltonian path. Obviously, a perimeter arc can be neither a
domino-axis nor a domino-black-end arc in case of any domino tiling. The
total number of white perimeter arcs is $p+q-2$. Since a Hamiltonian path
contains $pq-1$ arcs, there are exactly $pq-1-(p+q-2)=\allowbreak \left(
p-1\right) \left( q-1\right) $ arcs in any Hamiltonian path which are
neither white perimeter arcs, nor domino-axis arcs nor domino-black-end
arcs. Given a domino tiling, we call all arcs as \emph{domino-black-side}
arcs which are neither white perimeter arcs, nor domino-axis arcs nor
domino-black-end arcs.

Observe that around any black square, any Hamiltonian path contains exactly
two opposite arcs. Therefore from all black squares a Hamiltonian path
contains exactly $\left( p-1\right) \left( q-1\right) $ arcs.

Let $A$ be the set of all arcs in the grid which are not white perimeter
edges. We have that 
\[
\left \vert A\right \vert =2pq-q-p-(p+q-2)=2\left( p-1\right) \left(
q-1\right) 
\]%
(In Figure 6 we find an illustration for $p=4$, $q=3$.) Each Hamiltonian
path contains exactly half of the arcs in $A$.

Each domino tiling also takes exactly a quarter of the arcs of $A$ as
domino-black-end arcs and exactly a quarter of the arcs as domino-axis arcs.
And these three parts must be pairwise disjoint if and only if the
Hamiltonian path bisects neither domino of the tiling.

We obtained that a domino tiling determines the corresponding Hamiltonian
path. On the other hand, any Hamiltonian path determines for each black
square that in a domino tiling the domino containing the black square is in
vertical position or in horizontal position. This leaves one or two choices
for the white neighbor. However, if there are two choices, then both are
vertical or both are horizontal. In case of the bottom-left black square
there is only one choice. (In Figures 7.a and 7.b we find illustrations for $%
p=9$, $q=9$.)

Given a Hamiltonian path $H$ we make a bipartite graph $B_{H}$ as follows:
The black squares form one color class of the vertices in $B_{H}$ and the
white squares form the other color class of the vertices. A black and a
white square will be adjacent in $B_{H}$ if they are neighbors and the two
squares form such a domino whose position is allowed in the previous sense.
Observe that the Hamiltonian path allows no cycle in $B_{H}$. Therefore,
there exists at most one perfect matching in the bipartite graph. Therefore
a given Hamiltonian path $H$ allows at most one domino tiling whose dominoes
are all avoided by $H$. This completes the proof of Theorem 1.

\bigskip \noindent \textbf{References}

\smallskip \noindent [Ba09] Bang-Jensen, J.\ and Gutin, G.: \emph{Digraphs:
theory, algorithms and applications}, Second edition, Springer, 2009.

\smallskip \noindent [Gi07] Gindikin, S. (translated from Russian by Suchat,
A.): \emph{Tales of mathematicians and physicists}, Springer, 2007.

\smallskip \noindent [Ka61] Kasteleyn, P.\ W.: \emph{The statistics of
dimers on a lattice}, Physica \textbf{27} (1961) 1209--1225.

\smallskip \noindent [PL90] Prusinkiewicz, P.\ and Lindenmayer, A.: \emph{%
The algorithmic beauty of plants}, Springer, 1990.

\smallskip \noindent [PLF91] Prusinkiewicz, P., Lindenmayer, A., and
Fracchia, F. D.: \emph{Synthesis of space-filling curves on the square grid}%
. In: H.-O.\ Peitgen, J.\ M.\ Henriques \& L.\ F.\ Penedo, eds.: \emph{%
Fractals in the fundamental and applied sciences}, North-Holland (1991)
341--366.

\smallskip \noindent [Sa94] Sagan, H., \emph{Space-filling curves},
Springer, 1994.

\smallskip \noindent [Te61] Temperley, H.\ N.\ V.\ and Fisher, M.\ E.: \emph{%
Dimer problem in statistical mechanics--an exact result}, Phil.\ Mag.\ 
\textbf{68} (1961) 1061--1063.

\smallskip \noindent [We01] West, D.\ B.: \emph{Introduction to graph theory}%
, 2nd ed., Pearson Education, 2001.

\bigskip

\newpage \noindent \textbf{Figures} \bigskip

\[
\begin{picture}(180,130)
\put( 10  , 10  ){\line(1,0){160}}
\put( 10  , 90  ){\line(1,0){160}}
\put( 10  , 50  ){\line(1,0){160}}
\put( 10  ,130  ){\line(1,0){160}}
\put( 10  , 10  ){\line(0,1){120}}
\put( 50  , 10  ){\line(0,1){120}}
\put( 90  , 10  ){\line(0,1){120}}
\put(130  , 10  ){\line(0,1){120}}
\put(170  , 10  ){\line(0,1){120}}
\put( 10  , 10  ){\circle*{4}}
\put( 50  , 10  ){\circle*{4}}
\put( 90  , 10  ){\circle*{4}}
\put(130  , 10  ){\circle*{4}}
\put(170  , 10  ){\circle*{4}}
\put( 10  , 50  ){\circle*{4}}
\put( 50  , 50  ){\circle*{4}}
\put( 90  , 50  ){\circle*{4}} 
\put(130  , 50  ){\circle*{4}}
\put(170  , 50  ){\circle*{4}}
\put( 10  , 90  ){\circle*{4}}
\put( 50  , 90  ){\circle*{4}}
\put( 90  , 90  ){\circle*{4}} 
\put(130  , 90  ){\circle*{4}}
\put(170  , 90  ){\circle*{4}}
\put( 10  ,130  ){\circle*{4}}
\put( 50  ,130  ){\circle*{4}}
\put( 90  ,130  ){\circle*{4}}
\put(130  ,130  ){\circle*{4}}
\put(170  ,130  ){\circle*{4}}
\put( 27  ,  9.7){\makebox{$_>$}} 
\put( 67  ,  9.7){\makebox{$_>$}}  
\put(107  ,  9.7){\makebox{$_>$}}  
\put(147  ,  9.7){\makebox{$_>$}}  
\put( 26  , 49.7){\makebox{$_<$}} 
\put( 66  , 49.7){\makebox{$_<$}}  
\put(106  , 49.7){\makebox{$_<$}}  
\put(146  , 49.7){\makebox{$_<$}} 
\put( 27  , 89.7){\makebox{$_>$}} 
\put( 67  , 89.7){\makebox{$_>$}}  
\put(107  , 89.7){\makebox{$_>$}}  
\put(147  , 89.7){\makebox{$_>$}} 
\put( 26  ,129.7){\makebox{$_<$}} 
\put( 66  ,129.7){\makebox{$_<$}}  
\put(106  ,129.7){\makebox{$_<$}}  
\put(146  ,129.7){\makebox{$_<$}} 
\put(  7  ,  30  ){\makebox{$_\wedge$}} 
\put( 47  ,  25  ){\makebox{$_\vee$}} 
\put( 87  ,  30  ){\makebox{$_\wedge$}}
\put(127  ,  25  ){\makebox{$_\vee$}} 
\put(167  ,  30  ){\makebox{$_\wedge$}}
\put(  7  ,  70  ){\makebox{$_\wedge$}} 
\put( 47  ,  65  ){\makebox{$_\vee$}} 
\put( 87  ,  70  ){\makebox{$_\wedge$}}
\put(127  ,  65  ){\makebox{$_\vee$}} 
\put(167  ,  70  ){\makebox{$_\wedge$}}
\put(  7  , 110  ){\makebox{$_\wedge$}} 
\put( 47  , 105  ){\makebox{$_\vee$}} 
\put( 87  , 110  ){\makebox{$_\wedge$}}
\put(127  , 105  ){\makebox{$_\vee$}} 
\put(167  , 110  ){\makebox{$_\wedge$}}

\end{picture}
\]
\centerline{Figure 1.\ The directed graph $DGG_{5,4}$}

\bigskip

\[
\begin{picture}(180,130)
\put( 10  , 10  ){\line(1,0){ 80}}
\put(130  , 10  ){\line(1,0){ 40}}
\put( 10  , 10  ){\circle*{4}}
\put( 50  , 10  ){\circle*{4}}
\put( 90  , 10  ){\circle*{4}}
\put(130  , 10  ){\circle*{4}}
\put(170  , 10  ){\circle*{4}}
\put( 26  ,  9.7){\makebox{$_>$}} 
\put( 66  ,  9.7){\makebox{$_>$}} 
\put(146  ,  9.7){\makebox{$_>$}}   
\put( 10  , 90  ){\line(1,0){120}}
\put( 10  , 90  ){\circle*{4}}
\put( 50  , 90  ){\circle*{4}}
\put( 90  , 90  ){\circle*{4}} 
\put(130  , 90  ){\circle*{4}}
\put(170  , 90  ){\circle*{4}}
\put( 26  , 89.7){\makebox{$_>$}} 
\put( 66  , 89.7){\makebox{$_>$}}  
\put(106  , 89.7){\makebox{$_>$}}  
\put( 10  , 50  ){\line(1,0){ 80}}
\put( 10  , 50  ){\circle*{4}}
\put( 50  , 50  ){\circle*{4}}
\put( 90  , 50  ){\circle*{4}} 
\put(130  , 50  ){\circle*{4}}
\put(170  , 50  ){\circle*{4}}
\put( 23  , 49.7){\makebox{$_<$}} 
\put( 63  , 49.7){\makebox{$_<$}}  
\put( 10  ,130  ){\line(1,0){160}}
\put( 10  ,130  ){\circle*{4}}
\put( 50  ,130  ){\circle*{4}}
\put( 90  ,130  ){\circle*{4}}
\put(130  ,130  ){\circle*{4}}
\put(170  ,130  ){\circle*{4}}
\put( 23  ,129.7){\makebox{$_<$}} 
\put( 63  ,129.7){\makebox{$_<$}}  
\put(103  ,129.7){\makebox{$_<$}}  
\put(143  ,129.7){\makebox{$_<$}} 
\put( 10  , 50  ){\line(0,1){ 40}}
\put( 90  , 10  ){\line(0,1){ 40}}
\put(130  , 10  ){\line(0,1){ 80}}
\put(170  , 10  ){\line(0,1){120}}
\put( 87  , 30  ){\makebox{$_\wedge$}}
\put(127  , 25  ){\makebox{$_\vee$}} 
\put(167  , 30  ){\makebox{$_\wedge$}}
\put(  7  , 70  ){\makebox{$_\wedge$}} 
\put(127  , 65  ){\makebox{$_\vee$}} 
\put(167  , 70  ){\makebox{$_\wedge$}}
\put(167  ,110  ){\makebox{$_\wedge$}}
\end{picture}
\]
\centerline{Figure 2.a.\ A Hamiltonian path}

\bigskip

\[
\begin{picture}(180,130)
\put( 10  , 10  ){\line(1,0){160}}
\put( 10  , 50  ){\line(1,0){160}}
\put( 10  , 90  ){\line(1,0){160}}
\put( 10  ,130  ){\line(1,0){160}}
\put( 10  , 50  ){\line(0,1){ 40}}
\put(170  , 10  ){\line(0,1){ 40}}
\put(170  , 90  ){\line(0,1){ 40}}
\put( 10  , 10  ){\circle*{4}}
\put( 50  , 10  ){\circle*{4}}
\put( 90  , 10  ){\circle*{4}}
\put(130  , 10  ){\circle*{4}}
\put(170  , 10  ){\circle*{4}}
\put( 10  , 50  ){\circle*{4}}
\put( 50  , 50  ){\circle*{4}}
\put( 90  , 50  ){\circle*{4}} 
\put(130  , 50  ){\circle*{4}}
\put(170  , 50  ){\circle*{4}}
\put( 10  , 90  ){\circle*{4}}
\put( 50  , 90  ){\circle*{4}}
\put( 90  , 90  ){\circle*{4}} 
\put(130  , 90  ){\circle*{4}}
\put(170  , 90  ){\circle*{4}}
\put( 10  ,130  ){\circle*{4}}
\put( 50  ,130  ){\circle*{4}}
\put( 90  ,130  ){\circle*{4}}
\put(130  ,130  ){\circle*{4}}
\put(170  ,130  ){\circle*{4}}
\put( 27  ,  9.7){\makebox{$_>$}} 
\put( 67  ,  9.7){\makebox{$_>$}} 
\put(107  ,  9.7){\makebox{$_>$}}
\put(147  ,  9.7){\makebox{$_>$}}   
\put( 27  , 89.7){\makebox{$_>$}} 
\put( 67  , 89.7){\makebox{$_>$}}  
\put(107  , 89.7){\makebox{$_>$}} 
\put(147  , 89.7){\makebox{$_>$}}  
\put( 26  , 49.7){\makebox{$_<$}} 
\put( 66  , 49.7){\makebox{$_<$}}
\put(106  , 49.7){\makebox{$_<$}} 
\put(146  , 49.7){\makebox{$_<$}}   
\put( 26  ,129.7){\makebox{$_<$}} 
\put( 66  ,129.7){\makebox{$_<$}}  
\put(106  ,129.7){\makebox{$_<$}}  
\put(146  ,129.7){\makebox{$_<$}} 
\put(  7  , 70  ){\makebox{$_\wedge$}} 
\put(167  , 30  ){\makebox{$_\wedge$}}
\put(167  ,110  ){\makebox{$_\wedge$}}
\end{picture}
\]
\centerline{Figure 2.b.\ Another Hamiltonian path}

\bigskip

\[
\begin{picture}(180,130)
\put( 10  , 10  ){\line(1,0){160}}
\put(130  , 10  ){\line(1,0){ 40}}
\put( 10  , 90  ){\line(1,0){160}}
\put( 10  , 50  ){\line(1,0){ 80}}
\put( 10  ,130  ){\line(1,0){160}}
\put( 10  , 10  ){\line(0,1){120}}
\put( 90  , 10  ){\line(0,1){120}}
\put(130  , 10  ){\line(0,1){ 80}}
\put(170  , 10  ){\line(0,1){120}}
\put( 10  , 10  ){\circle*{4}}
\put( 50  , 10  ){\circle*{4}}
\put( 90  , 10  ){\circle*{4}}
\put(130  , 10  ){\circle*{4}}
\put(170  , 10  ){\circle*{4}}
\put( 10  , 90  ){\circle*{4}}
\put( 50  , 90  ){\circle*{4}}
\put( 90  , 90  ){\circle*{4}} 
\put(130  , 90  ){\circle*{4}}
\put(170  , 90  ){\circle*{4}}
\put( 10  , 50  ){\circle*{4}}
\put( 50  , 50  ){\circle*{4}}
\put( 90  , 50  ){\circle*{4}} 
\put(130  , 50  ){\circle*{4}}
\put(170  , 50  ){\circle*{4}}
\put( 10  ,130  ){\circle*{4}}
\put( 50  ,130  ){\circle*{4}}
\put( 90  ,130  ){\circle*{4}}
\put(130  ,130  ){\circle*{4}}
\put(170  ,130  ){\circle*{4}}
\end{picture}
\]
\centerline{Figure 3.a.\ A tiling corresponding to Fig.\ 2.a}

\bigskip

\[
\begin{picture}(180,130)
\put( 10  , 10  ){\line(1,0){160}}
\put( 10  , 50  ){\line(1,0){160}}
\put( 10  , 90  ){\line(1,0){160}}
\put( 10  ,130  ){\line(1,0){160}}
\put( 10  , 10  ){\line(0,1){120}}
\put( 90  , 10  ){\line(0,1){120}}
\put(170  , 10  ){\line(0,1){120}}
\put( 10  , 10  ){\circle*{4}}
\put( 50  , 10  ){\circle*{4}}
\put( 90  , 10  ){\circle*{4}}
\put(130  , 10  ){\circle*{4}}
\put(170  , 10  ){\circle*{4}}
\put( 10  , 50  ){\circle*{4}}
\put( 50  , 50  ){\circle*{4}}
\put( 90  , 50  ){\circle*{4}} 
\put(130  , 50  ){\circle*{4}}
\put(170  , 50  ){\circle*{4}}
\put( 10  , 90  ){\circle*{4}}
\put( 50  , 90  ){\circle*{4}}
\put( 90  , 90  ){\circle*{4}} 
\put(130  , 90  ){\circle*{4}}
\put(170  , 90  ){\circle*{4}}
\put( 10  ,130  ){\circle*{4}}
\put( 50  ,130  ){\circle*{4}}
\put( 90  ,130  ){\circle*{4}}
\put(130  ,130  ){\circle*{4}}
\put(170  ,130  ){\circle*{4}}
\end{picture}
\]
\centerline{Figure 3.b.\ A tiling corresponding to Fig.\ 2.b}

\bigskip
\[
\begin{picture}(180,130) 
\put( 10  , 10  ){\line(1,0){160}}
\put( 10  , 10  ){\circle*{4}}
\put( 50  , 10  ){\circle*{4}}
\put( 90  , 10  ){\circle*{4}}
\put(130  , 10  ){\circle*{4}}
\put(170  , 10  ){\circle*{4}}
\put( 10  , 50  ){\line(1,0){ 80}}
\put( 10  , 50  ){\circle*{4}}
\put( 50  , 50  ){\circle*{4}}
\put( 90  , 50  ){\circle*{4}} 
\put(130  , 50  ){\circle*{4}}
\put(170  , 50  ){\circle*{4}}
\put( 10  , 90  ){\line(1,0){160}}
\put( 10  , 90  ){\circle*{4}}
\put( 50  , 90  ){\circle*{4}}
\put( 90  , 90  ){\circle*{4}} 
\put(130  , 90  ){\circle*{4}}
\put(170  , 90  ){\circle*{4}}
\put( 10  ,130  ){\line(1,0){160}}
\put( 10  ,130  ){\circle*{4}}
\put( 50  ,130  ){\circle*{4}}
\put( 90  ,130  ){\circle*{4}}
\put(130  ,130  ){\circle*{4}}
\put(170  ,130  ){\circle*{4}}
\put( 10  , 10  ){\line(0,1){120}}
\put( 90  , 10  ){\line(0,1){120}}
\put(130  , 10  ){\line(0,1){ 80}}
\put(170  , 10  ){\line(0,1){120}}

\put(107  ,  9.7){\makebox{$_>$}}   
\put(107  , 89.7){\makebox{$_>$}}  
\put( 87  , 30  ){\makebox{$_\wedge$}}
\put(127  , 25  ){\makebox{$_\vee$}} 
\put( 87  , 70  ){\makebox{$_\wedge$}}
\put(127  , 65  ){\makebox{$_\vee$}} 

\put( 92  , 15  ){\makebox{$_1$}}   
\put( 92  , 50  ){\makebox{$_2$}}
\put( 92  , 82  ){\makebox{$_3$}}  
\put(132  , 82  ){\makebox{$_4$}} 
\put(132  , 50  ){\makebox{$_5$}}
\put(132  , 15  ){\makebox{$_6$}}

\end{picture}
\]
\centerline{Figure 4.\ A canonical numbering around a domino}

\bigskip

\[
\begin{picture}(180,130)

\put( 10  , 10  ){\line(1,0){160}}
\put( 10  , 90  ){\line(1,0){160}}
\put( 10  , 50  ){\line(1,0){160}}
\put( 10  ,130  ){\line(1,0){160}}
\put( 10  , 10  ){\line(0,1){120}}
\put( 50  , 10  ){\line(0,1){120}}
\put( 90  , 10  ){\line(0,1){120}}
\put(130  , 10  ){\line(0,1){120}}
\put(170  , 10  ){\line(0,1){120}}

\put( 10  , 10  ){\circle*{4}}
\put( 50  , 10  ){\circle*{4}}
\put( 90  , 10  ){\circle*{4}}
\put(130  , 10  ){\circle*{4}}
\put(170  , 10  ){\circle*{4}}

\put( 10  , 50  ){\circle*{4}}
\put( 50  , 50  ){\circle*{4}}
\put( 90  , 50  ){\circle*{4}} 
\put(130  , 50  ){\circle*{4}}
\put(170  , 50  ){\circle*{4}}

\put( 10  , 90  ){\circle*{4}}
\put( 50  , 90  ){\circle*{4}}
\put( 90  , 90  ){\circle*{4}} 
\put(130  , 90  ){\circle*{4}}
\put(170  , 90  ){\circle*{4}}

\put( 10  ,130  ){\circle*{4}}
\put( 50  ,130  ){\circle*{4}}
\put( 90  ,130  ){\circle*{4}}
\put(130  ,130  ){\circle*{4}}
\put(170  ,130  ){\circle*{4}}

\put( 10  , 10  ){\line(1,1){ 80}}
\put( 90  , 10  ){\line(1,1){ 40}}
\put( 10  , 90  ){\line(1,1){ 40}}
\put( 90  , 90  ){\line(1,1){ 40}}

\put( 20  , 10  ){\line(1,1){ 30}}
\put(100  , 10  ){\line(1,1){ 30}}
\put( 20  , 90  ){\line(1,1){ 30}}
\put(100  , 90  ){\line(1,1){ 30}}
\put( 10  , 20  ){\line(1,1){ 30}}
\put( 90  , 20  ){\line(1,1){ 30}}
\put( 10  ,100  ){\line(1,1){ 30}}
\put( 90  ,100  ){\line(1,1){ 30}}

\put( 50  , 60  ){\line(1,1){ 30}}
\put(130  , 60  ){\line(1,1){ 30}}
\put( 60  , 50  ){\line(1,1){ 30}}
\put(140  , 50  ){\line(1,1){ 30}}

\put( 50  , 70  ){\line(1,1){ 20}}
\put(130  , 70  ){\line(1,1){ 20}}
\put( 30  , 90  ){\line(1,1){ 20}}
\put(110  , 90  ){\line(1,1){ 20}}
\put( 40  , 90  ){\line(1,1){ 10}}
\put(120  , 90  ){\line(1,1){ 10}}
\put( 30  , 10  ){\line(1,1){ 20}}
\put(110  , 10  ){\line(1,1){ 20}}
\put( 40  , 10  ){\line(1,1){ 10}}
\put(120  , 10  ){\line(1,1){ 10}}

\put( 90  ,110  ){\line(1,1){ 20}}
\put( 70  , 50  ){\line(1,1){ 20}}
\put( 80  , 50  ){\line(1,1){ 10}}
 
\put( 10  , 110 ){\line(1,1){ 20}}
\put( 10  ,  30 ){\line(1,1){ 20}}
\put( 90  ,  30 ){\line(1,1){ 20}}

\put( 10  ,  40 ){\line(1,1){ 10}}
\put( 90  ,  40 ){\line(1,1){ 10}}
\put( 10  , 120 ){\line(1,1){ 10}}
\put( 90  , 120 ){\line(1,1){ 10}}

\put( 50  ,  80 ){\line(1,1){ 10}}
\put(160  ,  50 ){\line(1,1){ 10}}
\put(130  ,  80 ){\line(1,1){ 10}}
\put(150  ,  50 ){\line(1,1){ 20}}
\put(130  ,  50 ){\line(1,1){ 40}}

\end{picture}
\]

\centerline{Figure 5.\ The grid as a chessboard}

\bigskip\bigskip

\[
\begin{picture}(180,130)

\put( 10  , 10  ){\line(1,0){ 40}}
\put( 90  , 10  ){\line(1,0){ 40}}
\put( 10  , 90  ){\line(1,0){160}}
\put( 10  , 50  ){\line(1,0){160}}
\put( 10  ,130  ){\line(1,0){ 40}}
\put( 90  ,130  ){\line(1,0){ 40}}
\put( 90  , 10  ){\line(0,1){ 40}}
\put( 10  , 10  ){\line(0,1){ 40}}
\put( 10  , 90  ){\line(0,1){ 40}}
\put( 50  , 10  ){\line(0,1){120}}
\put( 90  , 10  ){\line(0,1){120}}
\put(130  , 10  ){\line(0,1){120}}
\put(170  , 50  ){\line(0,1){ 40}}

\put( 10  , 10  ){\circle*{4}}
\put( 50  , 10  ){\circle*{4}}
\put( 90  , 10  ){\circle*{4}}
\put(130  , 10  ){\circle*{4}}
\put(170  , 10  ){\circle*{4}}

\put( 10  , 50  ){\circle*{4}}
\put( 50  , 50  ){\circle*{4}}
\put( 90  , 50  ){\circle*{4}} 
\put(130  , 50  ){\circle*{4}}
\put(170  , 50  ){\circle*{4}}

\put( 10  , 90  ){\circle*{4}}
\put( 50  , 90  ){\circle*{4}}
\put( 90  , 90  ){\circle*{4}} 
\put(130  , 90  ){\circle*{4}}
\put(170  , 90  ){\circle*{4}}

\put( 10  ,130  ){\circle*{4}}
\put( 50  ,130  ){\circle*{4}}
\put( 90  ,130  ){\circle*{4}}
\put(130  ,130  ){\circle*{4}}
\put(170  ,130  ){\circle*{4}}

\put( 10  , 10  ){\line(1,1){ 80}}
\put( 90  , 10  ){\line(1,1){ 40}}
\put( 10  , 90  ){\line(1,1){ 40}}
\put( 90  , 90  ){\line(1,1){ 40}}

\put( 20  , 10  ){\line(1,1){ 30}}
\put(100  , 10  ){\line(1,1){ 30}}
\put( 20  , 90  ){\line(1,1){ 30}}
\put(100  , 90  ){\line(1,1){ 30}}
\put( 10  , 20  ){\line(1,1){ 30}}
\put( 90  , 20  ){\line(1,1){ 30}}
\put( 10  ,100  ){\line(1,1){ 30}}
\put( 90  ,100  ){\line(1,1){ 30}}

\put( 50  , 60  ){\line(1,1){ 30}}
\put(130  , 60  ){\line(1,1){ 30}}
\put( 60  , 50  ){\line(1,1){ 30}}
\put(140  , 50  ){\line(1,1){ 30}}

\put( 50  , 70  ){\line(1,1){ 20}}
\put(130  , 70  ){\line(1,1){ 20}}
\put( 30  , 90  ){\line(1,1){ 20}}
\put(110  , 90  ){\line(1,1){ 20}}
\put( 40  , 90  ){\line(1,1){ 10}}
\put(120  , 90  ){\line(1,1){ 10}}
\put( 30  , 10  ){\line(1,1){ 20}}
\put(110  , 10  ){\line(1,1){ 20}}
\put( 40  , 10  ){\line(1,1){ 10}}
\put(120  , 10  ){\line(1,1){ 10}}

\put( 90  ,110  ){\line(1,1){ 20}}
\put( 70  , 50  ){\line(1,1){ 20}}
\put( 80  , 50  ){\line(1,1){ 10}}
 
\put( 10  , 110 ){\line(1,1){ 20}}
\put( 10  ,  30 ){\line(1,1){ 20}}
\put( 90  ,  30 ){\line(1,1){ 20}}

\put( 10  ,  40 ){\line(1,1){ 10}}
\put( 90  ,  40 ){\line(1,1){ 10}}
\put( 10  , 120 ){\line(1,1){ 10}}
\put( 90  , 120 ){\line(1,1){ 10}}

\put( 50  ,  80 ){\line(1,1){ 10}}
\put(160  ,  50 ){\line(1,1){ 10}}
\put(130  ,  80 ){\line(1,1){ 10}}
\put(150  ,  50 ){\line(1,1){ 20}}
\put(130  ,  50 ){\line(1,1){ 40}}

\end{picture}
\]
\centerline{Figure 6.\ Without the white perimeter arcs}

\newpage

\ \ \ \ \ \ \ \ \ \ \ \begin{figure}[here]
\includegraphics[width=0.8\textwidth]{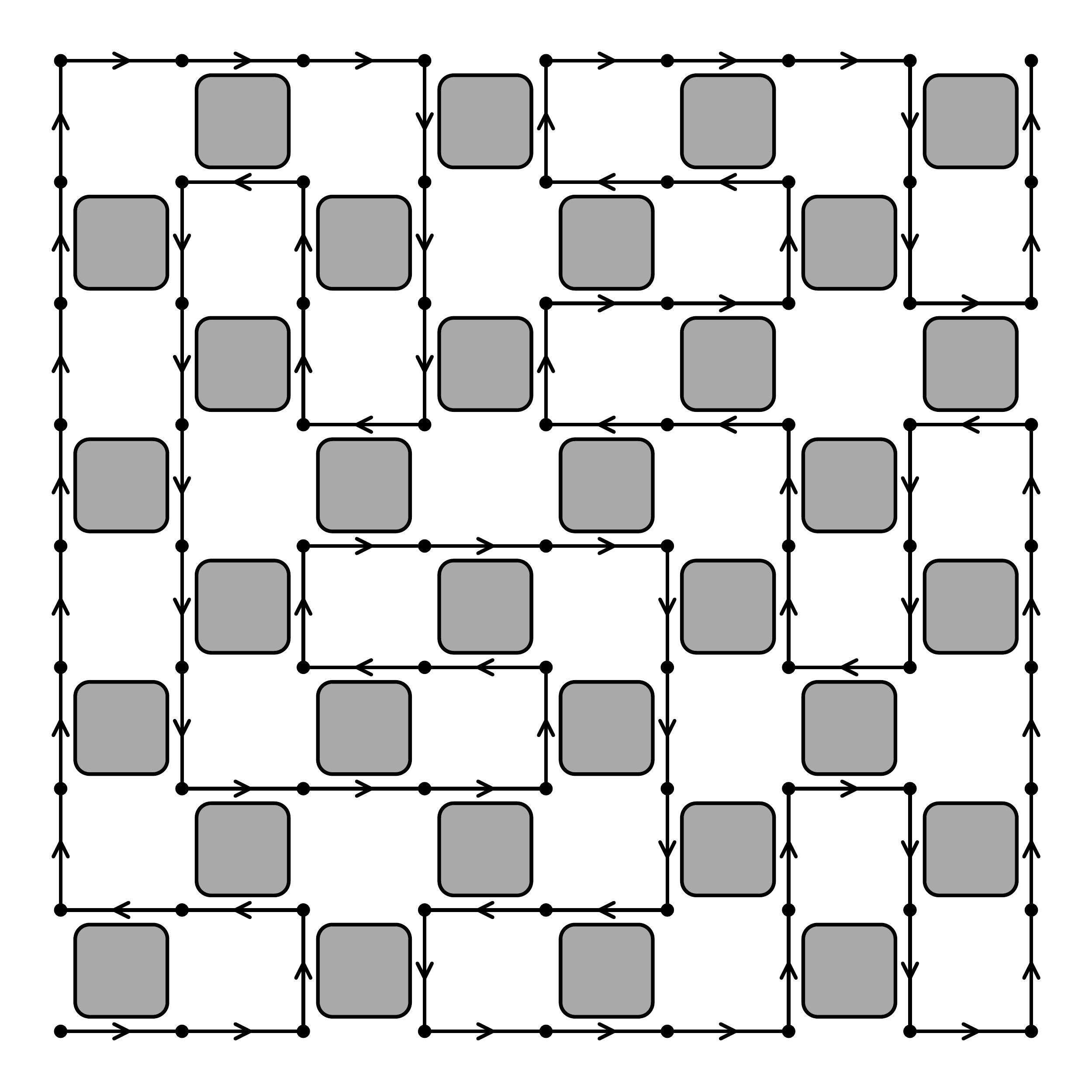}
\end{figure}

\centerline{Figure 7.a. The opposite arcs of the black squares and the white}

\centerline{perimeter arcs form a Hamiltonian path covering a Euclidean disk.}

\newpage

\ \ \ \ \ \ \ \ \ \ \ \ \ \ \ \begin{figure}[here]
\includegraphics[width=0.8\textwidth]{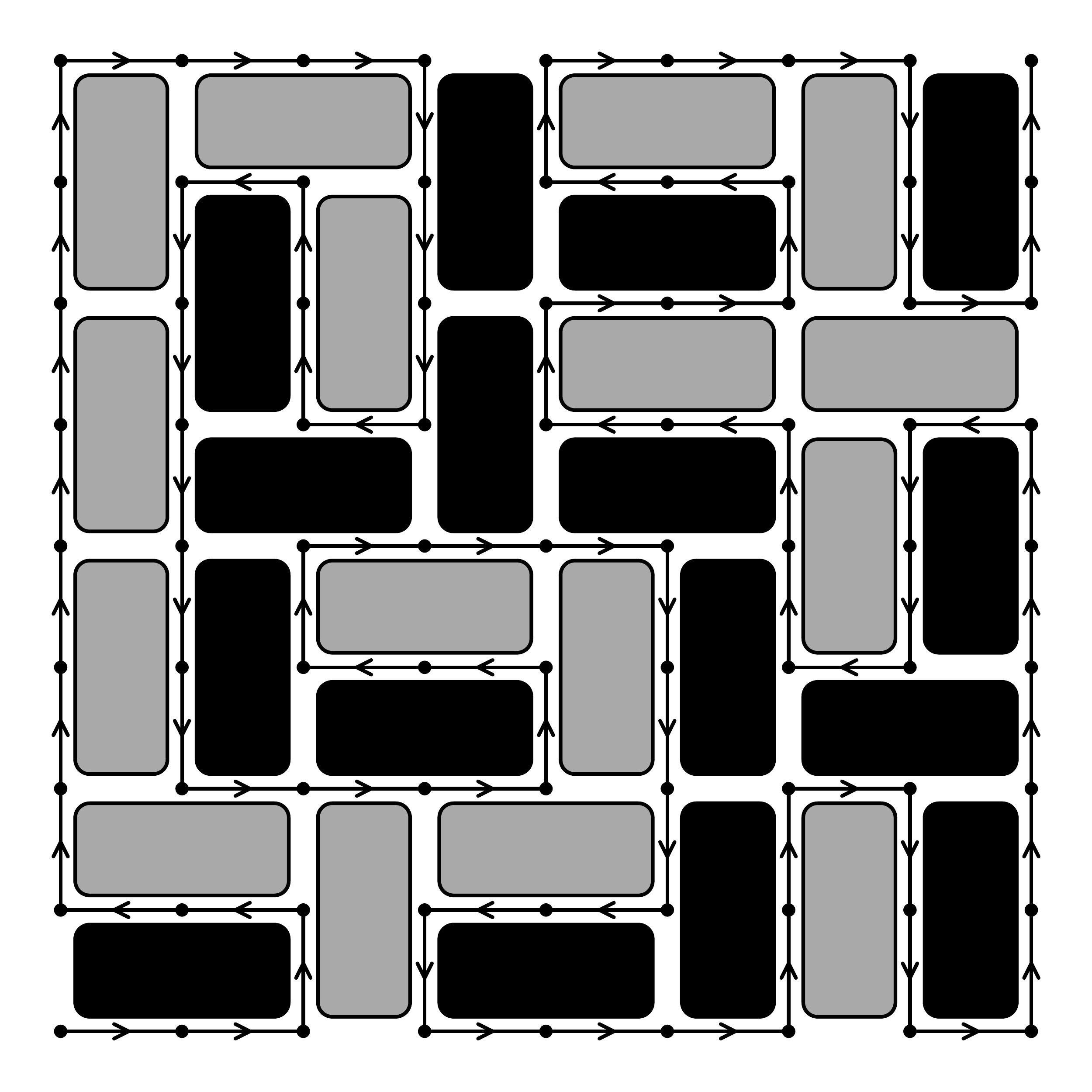}
\end{figure}

\centerline{Figure 7.b. A Hamiltonian path and the} 

\centerline{corresponding domino tiling.}

\newpage

\ \ \ \ \ \ \ \ \ \ \ \begin{figure}[here]
\includegraphics[width=0.7\textwidth]{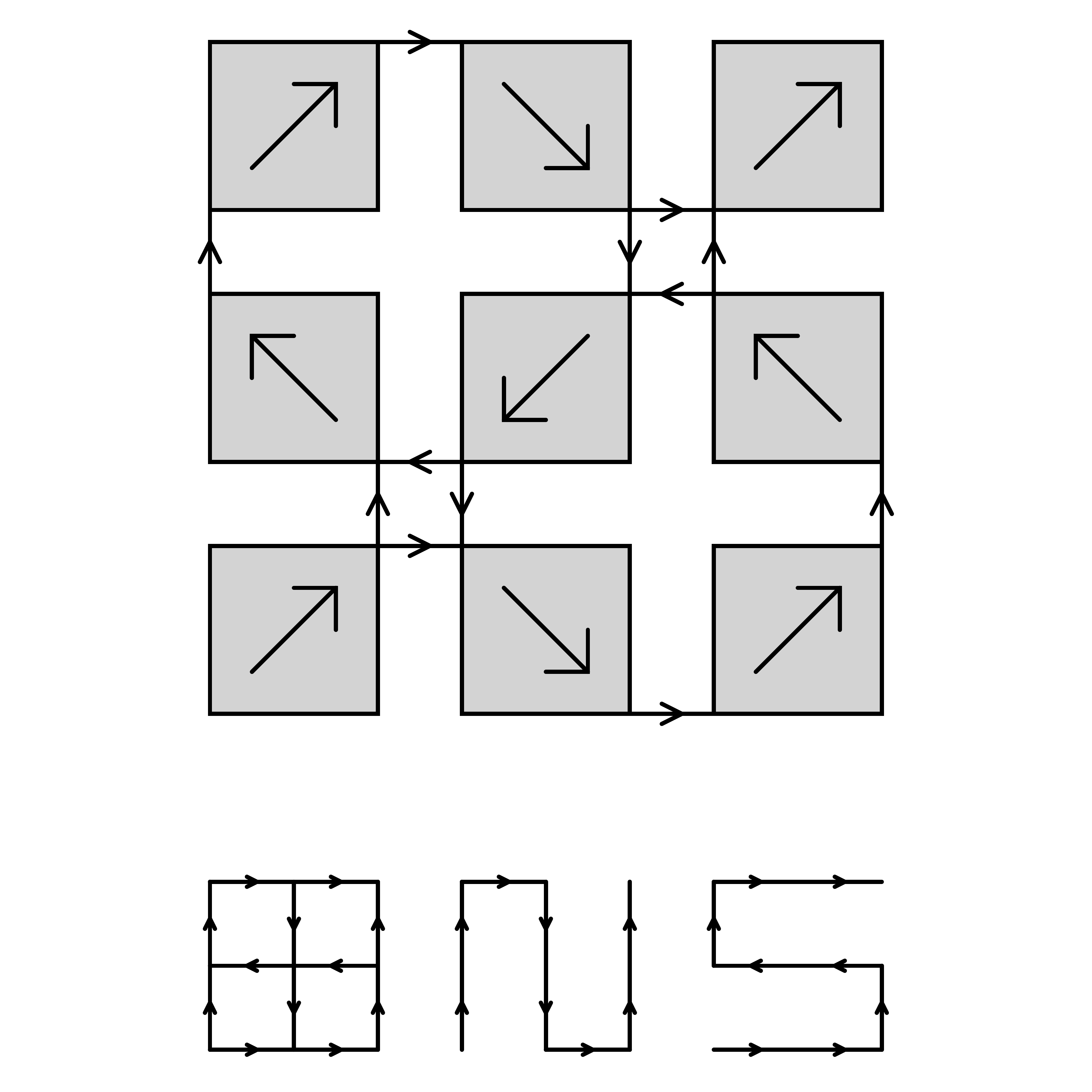}
\end{figure}

\centerline{Figure 8.a. Possible connecting edges between self similar parts}

\centerline{in ortho-tile type, node-rewriting FASS-curves.}

\centerline{(Grey squares grow out from the nodes of the path.)}

\medskip

\centerline{Figure 8.b. Directed Grid Graph: $DGG_{3,3}$}

\centerline{and the 2 possible Hamiltonian paths on it.}

\end{document}